\documentclass[a4paper,11pt]{extarticle}

\usepackage[left=3cm,right=1.5cm,top=2cm,bottom=2cm,bindingoffset=0cm]{geometry}
\usepackage{setspace} 

\usepackage{cmap}					
\usepackage[T2A]{fontenc}			
\usepackage[utf8]{inputenc}			
\usepackage[english]{babel}	

\usepackage{geometry} 
\usepackage{mathrsfs} 
\usepackage{graphicx}
\usepackage{amsmath,amsthm,amsfonts,amssymb} 
\usepackage{icomma} 

\usepackage{array,tabularx,tabulary,booktabs} 
\usepackage{longtable}  
\usepackage{multirow} 
\usepackage{hhline}

\usepackage{graphicx}  
\graphicspath{{images/}}  
\usepackage{wrapfig} 
\usepackage{epstopdf} 

\usepackage[usenames]{color}
\usepackage{color}
\usepackage{colortbl}

\usepackage{titlesec} 
\titleformat{\section}{\filcenter\normalfont\Large\bfseries}{\thesection.}{0.2em}{} 

\usepackage{indentfirst}

\begin{document}

\begin{center}
\textbf{REGULARIZATION OF THE DISCRETE SOURCE PROBLEM IN THE NONLINEAR DIFFUSIVE-LOGISTIC EQUATION
\footnote{The research was carried out within the state assignment of Ministry of Science and Higher Education of the Russian Federation (theme №~FENG-2023-0004, "Analytical and numerical study of inverse problems on recovering parameters of atmosphere or water pollution sources and (or) parameters of media").}}


\vspace{1em}
\textbf{Olga Krivorotko$^{1,2a}$, Tatiana Zvonareva$^{2b}$}


\vspace{1em}
\textit{
$^1$Yugra University, Khanty-Mansiysk, 628012 Russia,\\
$^2$Sobolev Institute of Mathematics SB RAS, Novosibirsk, 630090 Russia\\
e-mail: $^a$krivorotko.olya@mail.ru, $^b$t.zvonareva@g.nsu.ru}
\end{center}

A numerical algorithm for regularization of the solution of the source problem for the diffusion-logistic model based on information about the process at fixed moments of time of integral type has been developed. The peculiarity of the problem under study is the discrete formulation in space and impossibility to apply classical algorithms for its numerical solution. The regularization of the problem is based on the application of A.N.~Tikhonov's approach and a priori information about the source of the process. The problem was formulated in a variational formulation and solved by the global tensor optimization method. It is shown that in the case of noisy data regularization improves the accuracy of the reconstructed source.

\vspace{1em}
\textbf{Keywords:} diffusion-logistic model, social processes, inverse problem, optimization, regularization, a priori information.

\section*{INTRODUCTION}
The paper deals with the numerical solution of the source problem for the diffusion-logistic equation with inverse time
\begin{eqnarray}\label{eq:common}
    -\tilde{u}_t = a^2\Delta \tilde{u} + F(\tilde{u},x,t) + \varphi(x)\delta(t),
\end{eqnarray}
arising in the description of ecological, epidemiological and socio-economic processes with source function $\varphi(x)\in C^1 (\mathbb{R}^n)$ and nonlinear right-hand side $f\in C(t\ge 0)$~\cite{Pyatkov_2011, Viguerie_2020, Aristov_2021}.
Here $\delta (t)$ is the Dirac delta function characterizing the instantaneous propagation of the process under study. The problem consists in determination of the source function $\varphi(x)$ from additional information about the process at fixed times:
\begin{eqnarray}\label{eq:common_inv}
    \tilde{u}(x_i, t_k; \delta) = f_k^\delta,\quad k=1,\ldots, N_2.
\end{eqnarray}

In this paper, the problem~(\ref{eq:common})-(\ref{eq:common_inv}) will be investigated in the space of discrete values of the variable $x$, which can be responsible for the number of <<friendly ties>> from the source to the propagation of information in social processes, coordinates of pollution sources in ecological processes, etc. To formulate a continuous formulation of the problem, it is necessary to apply the interpolation of functions on the variable $x$, which is already ill-posed problem (its solution is non-unique)~\cite{Kabanikhin_2009}. Application of classical algorithms for numerical solution of the problem may lead to instability of the obtained solution.

In this paper, the solution of the source problem will be reduced to the solution of the minimization problem of the A.N.~Tikhonov functional with a constant regularization parameter. The solution of the minimization problem will be obtained by applying the global tensor optimization approach~\cite{Zheltkova_2018}. As an additional regularization, a priori information about a part of the solution (source) was used.

\section{INVERSE PROBLEM STATEMENT}

This paper investigates the numerical solution of a diffusion-logistic model~\cite{Wang_2013}:
\begin{eqnarray}\label{eq:DLM}
\left\{\begin{array}{ll}
    u_t = Du_{xx} + \left( 1 - \dfrac uK\right)r(t)u, & l_1\le x\le l_2, \, t\ge 1, \\
    u(x,1) = \varphi(x), & l_1\le x\le l_2, \\
    u_x (l_1,t) = u_x (l_2,t) = 0, & t\ge 1,
\end{array}\right.
\end{eqnarray}
obtained from the generalized equation~(\ref{eq:common}) in the case of one-dimensional space on the variable $x$, $D=a^2$, $F(u,x,t) = \left( 1 - \dfrac uK\right)r(t)u$~\cite{Vladimirov_1971}.

We consider the information dissemination process in online social networks, where $u(x,t)$ is the density of users involved in the information dissemination process, $D$ is the popularity of the information, $K$ is a constant describing the maximum throughput of the network, $r(t)$ is the growth rate of the number of active users, $\varphi(x)$~is the initial density function of the involved users. By distance $x$ we mean the minimum number of links from a source ($x = 0$) to an engaged
user. A more detailed description of the model is presented in~\cite{Zvonareva_2022}.

In the inverse problem, in addition to the function $u(x,t)$, the unknown is the initial density function $\varphi(x)$, in which the distance variable $x$ is a discrete quantity. The inverse problem is to determine the functions $u(x,t)$ and $\varphi(x)$ of the model~(\ref{eq:DLM}) from the noisy additional information about the process of integral-type:
\begin{equation}\label{eq:ip_data}
    f_k^{\delta} = f_k + \dfrac{\delta\cdot f_k\cdot\gamma_k}{100}, \quad f_{k} = \sum\limits_{i=1}^{N_{1}}u(x_{i},t_{k}), \quad k=1,\dots,N_{2},
\end{equation}
where $\delta$ is the noise level, $\gamma_k\sim\mathcal{N}(0,1)$ is the standard normal random variable. Data of the form~(\ref{eq:ip_data}) describe the density of users involved in the process of information dissemination at fixed moments of time $t_{k}$, $k=1,\dots,N_{2}$.

Note that the inverse problem~(\ref{eq:DLM})-(\ref{eq:ip_data}) is ill-posed~\cite{Kabanikhin_2009}, i.e., its solution may be non-unique and/or unstable.

Let $u(x,t)\in L^{2}((l_1,l_2)\times (1,+\infty))$. Then the inverse problem can be reduced to the corresponding minimization problem of the following regularizing functional A.~N.~Tikhonov~\cite{Tikhonov_1943, Yagola_1983}:
\begin{eqnarray}\label{eq:T}
    T(q) = \dfrac{T-1}{N_{2}}\sum_{k=1}^{N_{2}}\left|\sum_{i=1}^{N_{1}}u(x_{i},t_{k};q) - f_{k}^{\delta}\right|^2 + \alpha \sum\limits_{j = 0}^{d - 1}\left|q_j - q^0_j\right|^2,
\end{eqnarray}
where $u(x,t;q)$ is the solution of the direct problem for the initial density function $\varphi(x)$ defined from the set of parameters $q = (q_{0}, q_{1},\dots, q_{d - 1})$, $d=6$ is the number of parameters. Namely, $\varphi(x_i) = q_{i - 1}$, $i = 1,\dots, d$, and on each segment $[x_i, x_{i + 1}]$, $i = 1,\dots, d - 1$ the function $\varphi$ is a polynomial of degree three $\varphi_i(x)$.

In~\cite{Zvonareva_2023.1}, the solution to the inverse problem with noise-free data of the form~(\ref{eq:ip_data}) is obtained by applying the tensor optimization approach~\cite{Zheltkova_2018} to the functional~(\ref{eq:T}) in the absence of a regularizing term ($\alpha=0$). In this paper, the tensor optimization method will be applied to the minimization of the A.~N.~Tikhonov~(\ref{eq:T}) functional with a constant regularization parameter $\alpha$, and the use of a priori information about the solution of the inverse problem will be added.

\section{NUMERICAL EXPERIMENTS}

The direct problem~(\ref{eq:DLM}) is solved using an explicit finite-difference scheme with approximation order $O(\tau+h^2)$. A more detailed description of the numerical solution is presented in~\cite{Zvonareva_2022}.

In the A.~N.~Tikhonov functional~(\ref{eq:T}), the components of the vector $q^0_j\in [0,6]$ decrease linearly and are uniformly located on the segment $[0, 6]$, that is, $q^0 = (6, \, 4.8, \, 3.6, \, 2.4, \, 1.2, \, 0)$. 

We use the synthetic data $f_{k}$, that were taken from the solution of the direct problem at points $t = 3,\dots, 24$ and summed over $x = 1,\dots, 5$, i.e., $N_{1} = 5$ and $N_{2} = 22$. The values of the parameters $q_{ex} = (5.8, \, 1.7, \, 1.9, \, 1, \, 0.95, \, 0.7)$ emulate the real data of the news site Digg.com presented in~\cite{Wang_2013}.

The relative error for inverse problem solution was calculated using the formula:
\begin{align*}
    err_m = \dfrac{|q_{ex} - q^m|}{|q_{ex}|}.
\end{align*}
Here $q_{ex}$ are the exact values of the discretized function $\varphi(x)$, and $q^m$ is the approximation of the inverse problem solution, which corresponds to the minimum of the functional $T(q^m)$~(\ref{eq:T}).

The problem of minimizing the functional~(\ref{eq:T}) is solved using the tensor train optimization method (TT)~\cite{Zheltkova_2018, Zvonareva_2023.1}, the idea of which is to represent the functional $T(q^m)$ as a tensor on the parameter space. The algorithm of TT is as follows:
\begin{itemize}
   \item[0.] TT input: Lower and upper bounds of the solution space $b_{min}$ and $b_{max}$, number of parameters (dimensionality of the solution space) $d$, number of nodes along all directions $n$, maximum possible rank of cars $r_{max}$, number of iterations $N_{TT}$, initial shift of the functional $\alpha$, mapping function $h(J(q) - \alpha)$.
    \item[1.] Introduce a grid with $n$ nodes in each of the $d$ directions.
    \item[2.] Cycle from 1 execute to $d - 1$:
    \begin{itemize}
        \item[$\circ$] Using the grid values and the $\hat q_{i - 1}$ obtained in the previous step, generate $\hat q_i$.
    \end{itemize}
    \item[3.] As long as the number of iterations < $N_{TT}$, execute the
    \begin{itemize}
        \item[$\circ$] Cycle from 1 execute to $d - 1$:
        \begin{itemize}
            \item[$\bullet$] Based on $\hat q_{i - 1}$ and $\hat q_i$, generate a set of potential solutions $M$ and update the shift $\alpha$.
            \item[$\bullet$] Remember the best solution $q_{best}$.
            \item[$\bullet$] Represent the array of values of the function $h(\tilde q)$, $\tilde q\in M$, as a tensor.
            \item[$\bullet$] Compute the tensor approximation in TT-format.
            \item[$\bullet$] Using the grid values and $\hat q_{i - 1}$, generate $\hat q_i$.
        \end{itemize}
    \end{itemize}
\end{itemize}

\begin{figure}[!h]
    {\includegraphics[width=1\linewidth]{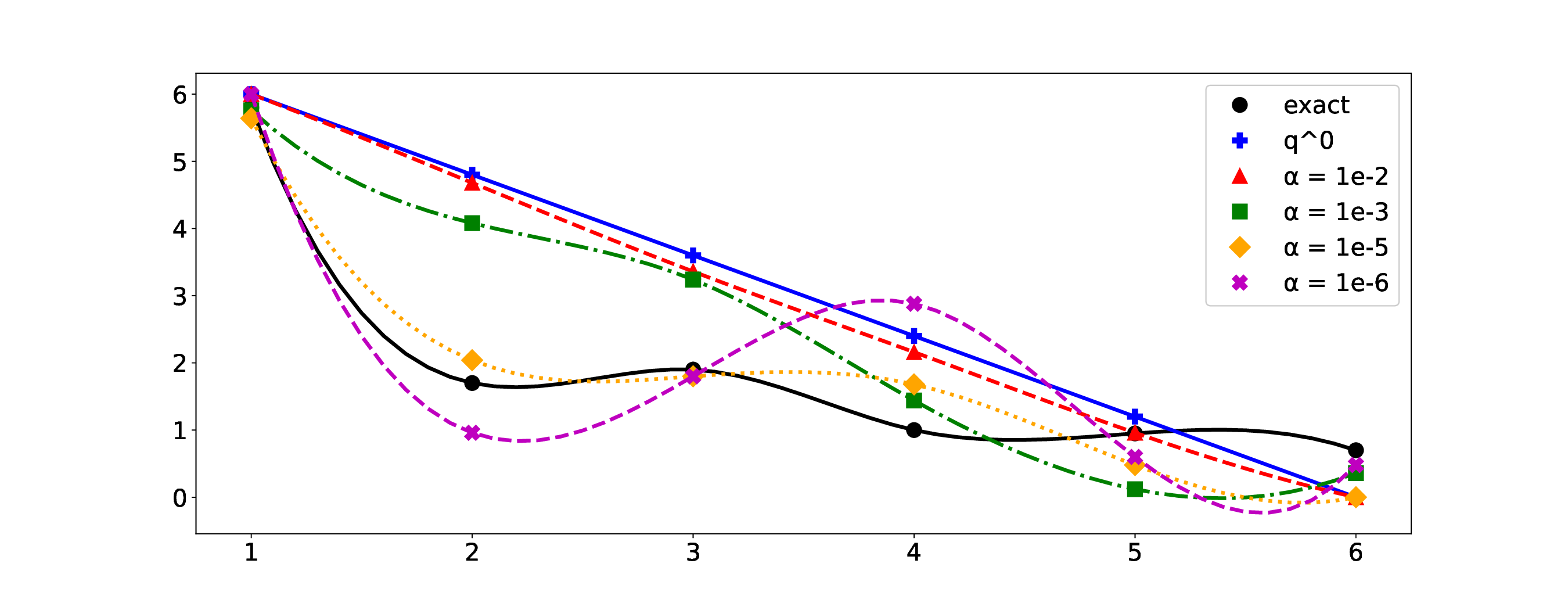}}
    \caption{The recovered initial density function of active users $\varphi(x)$ at $\delta = 10$\%. The continuous black line with circles represents the exact solution of the inverse problem, the blue line with crosses represents the initial approximation $q^0$.}
    \label{fig:al_del=10}
\end{figure}

Fig.~\ref{fig:al_del=10} shows the results of minimizing the $T(q)$~(\ref{eq:T}) functional by the TT method for the data $f_k^{\delta}$~(\ref{eq:ip_data}) with $\delta = 10$\% for different values of the regularization parameter $\alpha$.

The figure shows that the algorithm approximates the solution to $q^0$ for the sufficiently large values of the regularization parameter $\alpha$. This happens due to the too large contribution of the regularization term to the value of the functional. Therefore, as $\alpha$ decreases, the solution approaches the exact solution. However, when the regularization parameter is strongly reduced, the solution again becomes far away from the exact solution, indicating that the inverse problem~(\ref{eq:DLM})-(\ref{eq:ip_data}) is ill-posed. In this case, the smallest error value is achieved using the regularization parameter $\alpha = 10^{-5}$.

In the case of $\alpha = 0$ for different noise levels $\delta$, the results presented in Table~\ref{tabl:del_al=10-2}. That is, even at the minimum noise level $\delta = 1$, the TT method does not find a solution close to the exact one, while having the smallest value of the functional. And the smallest value of the error is achieved at $\delta = 5$\%.

\begin{table}[!ht]\centering
    \caption{Values of error $err$ and functional $T$ for different $\delta$ (in \%).}
    \label{tabl:del_al=10-2}
    \begin{tabular}{c|c|c} 
        \hline
        $\delta$ & $err$ & $T$\\
        \hline
        1 & 5.97 & $6.51\cdot10^{-6}$ \\
        2 & 4.89 & $2.83\cdot10^{-5}$ \\
        5 & 2.61 & $1.08\cdot10^{-4}$ \\
        8 & 14.78 & $3.62\cdot10^{-4}$ \\
        10 & 13.63 & $7.12\cdot10^{-4}$ \\
        \hline
    \end{tabular}
\end{table}

\begin{figure}[!h]
    \center{\includegraphics[width=1\linewidth]{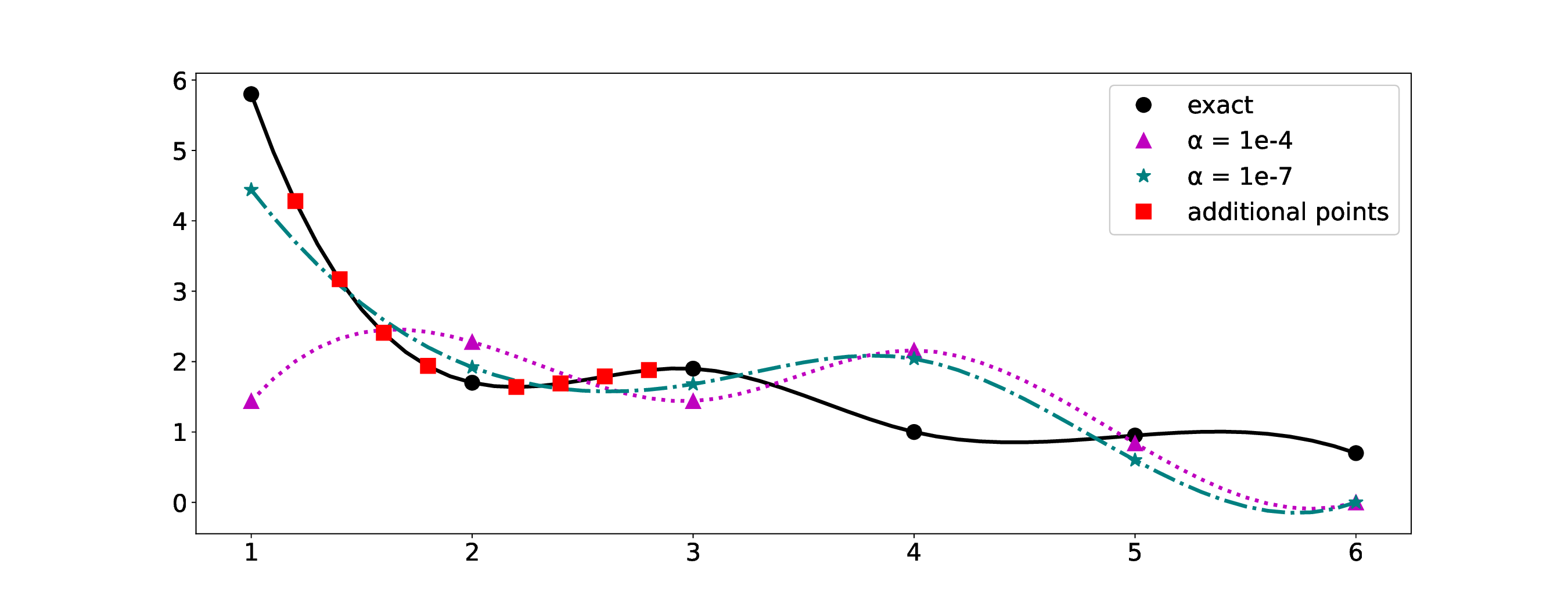}}
    \caption{Reconstructed initial density function of active users $\varphi(x)$ at $\delta = 10$ and $q^0 = 0$. The continuous black line with circles represents the exact solution of the inverse problem, the purple dashed line with triangles represents the solution obtained by the TT method with regularization parameter $\alpha = 10^{-4}$, the blue-green (teal) dashed line with stars represents the solution obtained by the TT method with regularization parameter $\alpha = 10^{-7}$, and the red squares represent additional points}.
    \label{fig:add_p}
\end{figure}

Similar to Fig.~\ref{fig:al_del=10} results were obtained for the case $q^0 = 0$. The TT method approximated the solution to zero at sufficiently large values of $\alpha$, but decreasing the value of the regularization parameter did not restore the behavior of the initial density function $\varphi(x)$ on the segment $x\in[1,2]$. Only the solution at $\alpha = 10^{-7}$ ''guesses'' the behavior of the function on this segment, and the value of the parameter is too small, which almost eliminates the influence of the regularizing parameter. Therefore, we decided to consider the inverse problem at $d = 14$ with the addition of a priori information at the distances $x\in [1, 3]$ closest to the propagation source, as presented in Fig.~\ref{fig:add_p}.

\begin{figure}[!ht]
    \center{\includegraphics[width=1\linewidth]{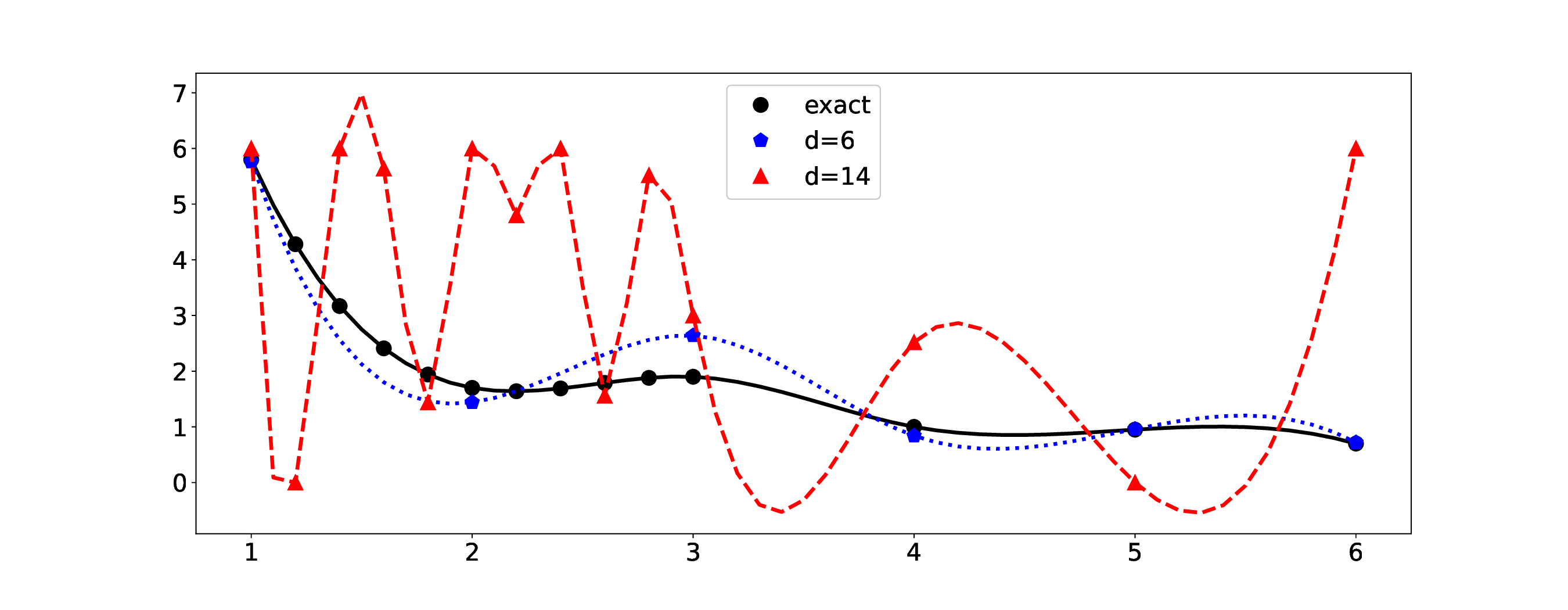}}\\ (a)
    \vspace{-1em}
    \center{\includegraphics[width=1\linewidth]{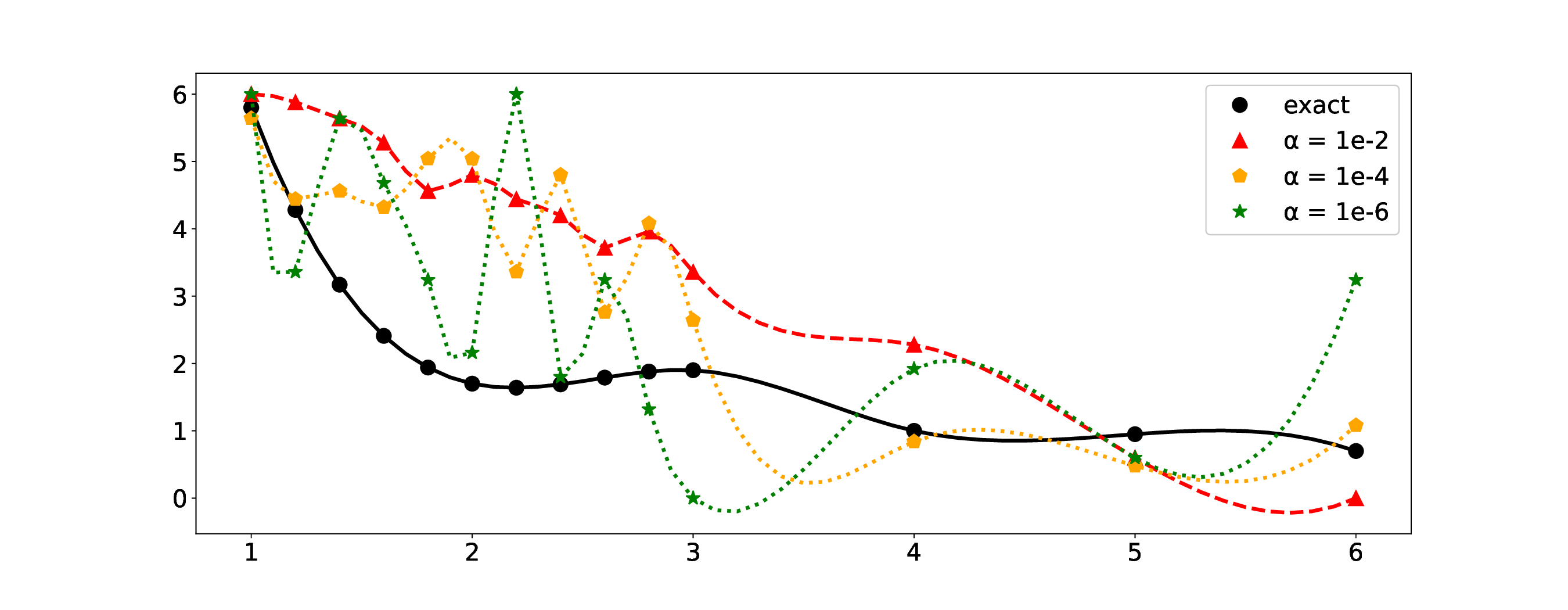}}\\ (b)
    \caption{Reconstructed initial density function of active users $\varphi(x)$ at $\delta = 0$ for the case $d = 14$. (a) comparison of the reconstruction results for the cases $d = 6$ and $d = 14$ at $\alpha = 0$; (b) reconstruction results of the source function for different values of the regularization parameter $\alpha$.} 
    \label{fig:d=14}
\end{figure}

Thus, we considered the problem of minimizing the A.~N.~Tikhonov~(\ref{eq:T}) functional at $\delta = 0$ for the case $d = 14$, $q_{ex} = (5.8, \, 4.28, \, 3.17, \, 2.41, \, 1.94, \, 1.7, \, 1.64, \, 1.69, \, 1.79, \, 1.88, \, 1.9, \, 1$, $0.95, \, 0.7)$. The results are shown in Fig.~\ref{fig:d=14}: (a) -- comparison of the reconstruction results for the cases $d = 6$ and $d = 14$ at $\alpha = 0$; (b) -- results of the source function reconstruction for different values of the regularization parameter $\alpha$. In Fig.~\ref{fig:d=14}(a) it can be seen that for $d = 6$ the solution is recovered quite close to exact, while for the case $d = 14$ there is not enough data (note that the initial density is recovered from $N_2=22$ values) to bring the solution close to exact. Applying regularization (Fig.~\ref{fig:d=14}(b) at $\alpha\not=0$) does not improve the solution. The above trend is observed, in which for sufficiently large values of $\alpha$ the solution approaches $q^0$, but decreasing the regularization parameter does not improve the solution.

\section*{CONCLUSION}

In this paper we study the regularization algorithm for the solution of the source recovery problem for the diffusion-logistic equation using the density of influenced users in dissemination of information in online social networks in time points (inverse problem). The initial formulation of the inverse problem is discrete in space, which leads to nonuniqueness and instability of the solution obtained by classical approaches and to the necessity of regularization. The inverse problem was reduced to the corresponding problem of minimization of the regularizing functional A.~N.~Tikhonov and numerically solved using the method of tensor optimization. The formulation of the problem of source determination by additional information of integral type with noisy data is considered. 

It is shown that the algorithm approximates the solution to the given initial approximation $q^0$ for the sufficiently large values of the regularization parameter $\alpha$. This occurs due to the too large contribution of the regularization term to the value of the functional. Therefore, when $\alpha$ decreases, the solution approaches the exact solution, but when the regularization parameter is strongly decreased, the solution is again distant from the exact solution, which indicates that the inverse problem is incorrect. In the case of $\alpha = 0$, even at the minimum noise level $\delta = 1$, the TT method does not find a solution close to the exact one.

The problem of recovering the vector of 14 parameters was also considered. It is shown that for 6 parameters the solution is recovered quite close to the exact one, while for the case of 14 parameters there is not enough data to bring the solution close to the exact one. And application of a priori information does not improve the solution sufficiently. In the future, we plan to use a combination of global optimization methods and gradient method, in which the descent parameter is matched with the regularization parameter.

\end{document}